\documentclass[12pt]{article}
\usepackage{e-jc}

\usepackage{amsmath,amsfonts}
\usepackage{asymptote}
\usepackage{cprotect}
\usepackage{enumerate}
\usepackage{graphicx}
\usepackage{mathtools}
\usepackage[capitalise,nameinlink]{cleveref}

\DeclarePairedDelimiter{\pa}{\lparen}{\rparen}

\begin{asydef}
int[] parse(string s) {
    int[] p = new int[length(s)];
    for (int i = 0; i < length(s); ++i) {
        int c = ascii(substr(s, i, 1));
        if (ascii("A") <= c && c <= ascii("Z"))
            p[i] = c - ascii("A") + 10;
        else if (ascii("0") <= c && c <= ascii("9"))
            p[i] = c - ascii("0");
        else
            p[i] = 0;
    }
    return p;
}
void drawperm(string s, real offset=0.0) {
    int[] p = parse(s);
    int len = p.length;
    for (int i = 0; i < len; ++i) {
        real lx = offset + 1.0 * i / len;
        real ry = 1.0 * p[i] / len;
        real d = 1.0 / len;
        fill((lx, ry)--(lx, ry - d)--(lx + d, ry - d)--(lx + d, ry)--cycle, gray);
    }
    draw((offset, 0)--(offset, 1)--(offset + 1, 1)--(offset + 1, 0)--cycle);
}
\end{asydef}

\dateline{Aug 18, 2020}{Jan 10, 2021}{TBD}
\Copyright{The authors. Released under the CC BY-ND license (International 4.0).}

\title{Limit Densities of Patterns in Permutation Inflations}
\author{
Tanya Khovanova\\
\small Department of Mathematics\\[-0.8ex]
\small Massachusetts Institute of Technology\\[-0.8ex]
\small Cambridge, MA, USA\\
\small\tt tanya@math.mit.edu\\
\and
Eric Zhang\\
\small Department of Mathematics\\[-0.8ex]
\small Harvard University\\[-0.8ex]
\small Cambridge, MA, USA\\
\small\tt ekzhang@college.harvard.edu}

\begin{document}

\maketitle

\begin{abstract}
Call a permutation $k$-inflatable if the sequence of its tensor products with uniform random permutations of increasing lengths has uniform $k$-point pattern densities. Previous work has shown that nontrivial $k$-inflatable permutations do not exist for $k \geq 4$. In this paper, we derive a general formula for the limit densities of patterns in the sequence of tensor products of a fixed permutation with each permutation from a convergent sequence. By applying this result, we completely characterize $3$-inflatable permutations and find explicit examples of $3$-inflatable permutations with various lengths, including the shortest examples with length $17$.
\end{abstract}

\section{Introduction}

In a broad sense, an object is called \emph{quasirandom} if, asymptotically, it has similar properties to random objects of the class it belongs to. The notion of quasirandomness has been studied for a variety of objects, including groups \cite{gowers07}, graphs \cite{chung89, lovasz08}, and $k$-uniform hypergraphs \cite{gowers06}. For permutations in particular \cite{cooper04, cooper08}, several different definitions of randomness are equivalent to a single concept of a \emph{quasirandom} permutation sequence. One such definition states that in the limit, as the lengths of the permutations in a sequence grow toward infinity, the densities of all pattern permutations of length $k$ approach $1/k!$.

The study of regularity and quasirandomness is an active area of research, and there have been fruitful prior results, such as upper bounds on packing density in \cite{sliacan18}. It was shown in \cite{kral12} that if a permutation sequence has uniform asymptotic densities for all $4$-patterns, then it is quasirandom. With this fact in mind, it is natural to ask whether there exist nontrivial permutation sequences that have uniform densities of $3$-patterns.

One method of constructing such permutation sequences is called \emph{random inflation}. We define the \emph{tensor product} of two permutations $\tau$ and $\gamma$ with lengths $n$ and $m$ to be a new permutation with length equal to the product $nm$. For each index $i \in [n]$ (the range of integers from $1$ to $n$) in the first permutation $\tau$, we replace it with $m$ numbers $n(\tau(i) - 1) + \gamma(j)$ for each $j \in [m]$. In an intuitive sense, the tensor product ``inflates'' each index of the first permutation by substituting for it $m$ indices in the order of the second permutation.

We can use the tensor product to build convergent permutation sequences. The main idea is to start from a base permutation and take its tensor product with a sequence of random permutations with lengths tending to infinity. The $n$-th term of this sequence is the tensor product of the base permutation with a permutation chosen uniformly at random from all elements of $S_n$.

A permutation is called \emph{$k$-inflatable} when this construction results in a convergent permutation sequence that has uniform densities of all $k$-patterns. This topic has been studied in the past. It was proved in \cite{kral12} that there are no $k$-inflatable permutations of length greater than $1$, for all $k \geq 4$. Also, previous work in \cite{cooper08} has analyzed $3$-inflatable and $4$-inflatable permutations.

We decided to write this paper when we discovered through experimentation that the example of a $3$-inflatable permutation of length 9 in \cite{cooper08} is wrong. With their methods as motivation, our primary contribution is a general method of computing asymptotic pattern densities in the random inflations of permutations. We used this formula to find conditions on the possible lengths of such permutations.

Combining these theoretical results with an extensive computational search, we correct the previous work by showing that the shortest length of a $3$-inflatable permutation is $17$, and we provide a non-exhaustive list of length-$17$ permutations that are $3$-inflatable. We also describe a method of generating an infinite number of $3$-inflatable permutations with increasing lengths. All of the code and accompanying files referenced in this paper are available online at \url{https://github.com/ekzhang/inflatable}.

We start with basic definitions in \cref{sec:preliminaries}, following the notation in \cite{albert05}. This includes formal definitions for permutation density and the convergence of permutation sequences, as well as a definition of quasirandomness for permutations.

In \cref{sec:density-tensor-product}, we derive a method for computing the limit density of a $k$-permutation $\pi$ in the tensor product of a permutation $\tau$ with permutations of a convergent sequence $\{ \gamma_j \}$, given the limit densities of patterns with lengths up to $k$ in that sequence.

In \cref{sec:density-inflation}, we examine the properties of random inflations, which are sequences given by the tensor product of a fixed permutation with uniform permutations of increasing length. By applying the theorem from the previous section, we can efficiently compute the limit densities of $k$-patterns in random inflations.

In \cref{sec:ksymmetry} we analyze $3$-point densities in random inflations. We show that $3$-inflatable permutations must have lengths in the set $\{0, 1, 17, 64, 80, 81\}$ modulo $144$. By computational search, we find examples of $3$-inflatable permutations for each of these residue classes of lengths, including the shortest instances with length $17$. We also discuss how rotational symmetry aids this computational search.

In \cref{sec:structure3} we study the structure of $3$-inflatable permutations. We show that if two permutations are both $3$-inflatable, then their tensor product is also $3$-inflatable.

\section{Preliminaries}
\label{sec:preliminaries}

A \textit{permutation} is an ordering of the elements of a set. We represent a permutation of length $n$ as an $n$-tuple of distinct positive integers, representing the image of $[n]$ under application of the permutation. For example, $(1, 2, 3)$ is the identity permutation of length 3, while $\tau = (3, 2, 1)$ is the permutation that swaps 1 and 3. For brevity, we usually leave out the parentheses. In the example above, we would write $\tau = 321$. The length of a permutation $\tau$ is denoted $|\tau|$, and the set of all permutations of length $n$ is $S_n$.

We are interested in patterns that are formed by subsets in a given permutation. Suppose we have an ordered $k$-tuple $(a_1, \ldots, a_k)$ of distinct positive integers. We say that this tuple is \textit{order-isomorphic} to a permutation $\pi$ if $a_i < a_j \iff \pi(i) < \pi(j)$ for all $i, j$.

\begin{definition}[Permutation density]
The \textit{density} of a $k$-pattern $\pi$ in a permutation $\tau$ is defined as the probability that a randomly-selected $k$-point subset of $\tau$ is order-isomorphic to $\pi$. We denote the density of $\pi$ in $\tau$ by $t(\pi, \tau)$. For example, $t(12, 132) = 2/3$.
\end{definition}

\begin{definition}[Convergent permutation sequence]
A sequence of permutations $\{ \tau_j \}$ is called \textit{convergent} if as $j$ increases, the length of the permutations $|\tau_j|$ approaches infinity, and for any permutation $\pi$, the densities $t(\pi, \tau_j)$ converge as a sequence of reals indexed by $j$. The limit \[
    \lim_{j  \to \infty} t(\pi, \tau_j)
\]
is called the \textit{limit density} of $\pi$ in the sequence $\{\tau_j\}$.
\end{definition}

This notion of density can be used to define the property of \emph{quasirandomness} for convergent permutation sequences. Notice that given any pattern permutation $\pi$, if one selects a permutation $\lambda$ uniformly at random from $S_n$ for $n \ge |\pi|$, we have by symmetry that $\mathbf E[t(\pi, \lambda)] = \frac{1}{|\pi|!}$, where $\mathbf E$ denotes expectation over choice of the random permutation $\lambda$. This motivates the following.

\begin{definition}[Quasirandom permutation sequence]
\label[definition]{def:quasirandom}
A convergent sequence of permutations $\{\tau_j\}$ is called \emph{quasirandom} if for every permutation $\pi$, \[
    \lim_{j \to \infty} t(\pi, \tau_j) = \frac{1}{|\pi|!}.
\]
\end{definition}

\section{Limit Densities in Permutation Products}
\label{sec:density-tensor-product}

We prove our first main result, which is a general method of computing the densities of patterns in the sequence of tensor products of a fixed permutation with a convergent sequence. However, before we can do this, we need some standard notation.

\subsection{Permutation Products}

These definitions were first introduced in \cite{albert05}. Some of our notation is slightly different (particularly for block decomposition), so we restate the definitions below.

\begin{figure}
    \centering
    \begin{asy}
    unitsize(3.5cm);
    drawperm("235461");
    draw((0,1/6)--(0,1/2)--(1/3,1/2)--(1/3,1/6)--cycle, blue + dotted + linewidth(1));
    draw((5/6,1)--(1/3,1)--(1/3,1/2)--(5/6,1/2)--cycle, blue + dotted + linewidth(1));
    draw((1,0)--(5/6,0)--(5/6,1/6)--(1,1/6)--cycle, blue + dotted + linewidth(1));
    label("$231[12, 213, 1] = 235461$", (0.5, -0.03), S);
    drawperm("312645", 1.2);
    draw(shift(1.2,0) * ((0,0)--(0,1/2)--(1/2,1/2)--(1/2,0)--cycle), blue + dotted + linewidth(1));
    draw(shift(1.2,0) * ((1,1)--(1,1/2)--(1/2,1/2)--(1/2,1)--cycle), blue + dotted + linewidth(1));
    label("$12[312] = 312645$", (0.5 + 1.2, -0.03), S);
    drawperm("561234", 1.2 * 2);
    draw(shift(1.2 * 2,0) * ((0,1)--(0,2/3)--(1/3,2/3)--(1/3,1)--cycle), blue + dotted + linewidth(1));
    draw(shift(1.2 * 2,0) * ((1/3,0)--(1/3,1/3)--(2/3,1/3)--(2/3,0)--cycle), blue + dotted + linewidth(1));
    draw(shift(1.2 * 2,0) * ((2/3,1/3)--(2/3,2/3)--(1,2/3)--(1,1/3)--cycle), blue + dotted + linewidth(1));
    label("$312[12] = 561234$", (0.5 + 1.2 * 2, -0.03), S);
    \end{asy}
    \caption{Examples of permutation inflation and tensor product.}
    \label{fig:inflation}
\end{figure}

\begin{definition}[Inflation]
Given a permutation $\tau$ of length $n$, the \textit{inflation} of $\tau$ with a sequence of $n$ permutations $\gamma_1, \ldots, \gamma_n$ is a permutation $\tau'$ of length $|\gamma_1| + \cdots + |\gamma_n|$ that consists of $n$ blocks, such that the $i$-th block is order-isomorphic to $\gamma_i$, and any restriction of $\tau'$ to one element in each block is order-isomorphic to $\tau$. This is denoted $\tau[\gamma_1, \ldots, \gamma_n]$.
\end{definition}

One example of inflation is $231[12, 213, 1] = 23\ 546\ 1$, where we add spaces for clarity. This is shown in \cref{fig:inflation}. Note how each element in the original permutation $\tau$ corresponds to a block of elements in the inflation.

There is a particular case of inflation that is often used. Rather than specifying the inflation of a permutation $\tau$ with sequence of $|\tau|$ different permutations $\gamma_1, \ldots, \gamma_\ell$, we have a special notation for the case when all of these blocks are the same, i.e., $\gamma_1 = \cdots = \gamma_\ell$.

\begin{definition}[Tensor product]
Given two permutations $\tau \in S_n$ and $\gamma \in S_m$, their \emph{tensor product} $\tau[\gamma]$ is a permutation of length $mn$ that consists of $n$ block of length $m$, where each block is order-isomorphic to $\gamma$, and the restriction of $\tau[\gamma]$ to one element in each block is order-isomorphic to $\tau$.
\end{definition}

For example, if $\tau = 12$ and $\gamma = 312$, then $\tau[\gamma] = 312\ 645$, shown in \cref{fig:inflation}. As another example, $\gamma[\tau] = 56\ 12\ 34$. Given a pattern $\pi$ and two permutations $\tau$, $\gamma$, it is possible to determine the densities of $\pi$ in $\tau[\gamma]$ based on densities of smaller patterns in $\tau$ and $\gamma$. We will use this to construct examples of permutation sequences with specified densities.

\begin{definition}[Block decomposition]
Given a permutation $\pi$, a \emph{block decomposition} of $\pi$ with order $\ell$ is a pair $(\sigma, b)$, such that $\sigma$ is a permutation of length $\ell$, and $b$ is a sequence of permutations $(b_1, \ldots, b_{\ell})$, such that $\pi = \sigma[b_1, \ldots, b_\ell]$.
\end{definition}

Finally, given any permutation $\pi$, we denote by $B(\pi)$ the set of all block decompositions of $\pi$. For example, if $\pi = 1243$, then $\pi$ has $6$ block decompositions: \begin{align*}
    B(1243) = \{
        &(1243, (1, 1, 1, 1)), \\
        &(132, (12, 1, 1)), \\
        &(123, (1, 1, 21)), \\
        &(12, (1, 132)), \\
        &(12, (12, 21)), \\
        &(1, (1243))
    \}
\end{align*}

\subsection{Convergent Sequences of Products}

Using block decompositions, we derive an explicit formula for the limit densities of patterns in the sequence of tensor products of a fixed permutation with each permutation from the convergent sequence $\{ \gamma_j \}_{j=1}^\infty$.

\begin{theorem}
\label{thm:tensor-limit-density}
For any permutation $\tau$ and convergent permutation sequence $\gamma_1, \gamma_2, \ldots$, the sequence of tensor products $\tau[\gamma_1], \tau[\gamma_2], \ldots$ is also convergent. Furthermore, the limit density of a pattern $\pi$ in this sequence is
\[
    \lim_{j\to\infty} t(\pi, \tau[\gamma_j]) = \frac{|\pi|!}{|\tau|^{|\pi|}} \sum_{(\sigma, b) \in B(\pi)} \pa*{ \binom{|\tau|}{|\sigma|} t(\sigma, \tau) \cdot \prod_{i=1}^{|\sigma|} \frac{\lim_{j\to\infty} t(b_i, \gamma_j)}{|b_i|!} }.
\]
\end{theorem}
\begin{proof}
The density of $\pi$ within $\tau[\gamma_j]$ equals the probability that when selecting $|\pi|$ independent random indices uniformly from $\tau[\gamma_j]$, the restriction of $\tau[\gamma_j]$ to these indices is order-isomorphic to $\pi$, conditioned on all indices being distinct. Note that as $j \to \infty$, the probability that all indices are distinct clearly approaches $1$. Therefore, \[
    \lim_{j \to \infty} t(\pi, \tau[\gamma_j]) = \lim_{j \to \infty} \frac{\#(\{|\pi|\text{ distinct indices order-isomorphic to } \pi\})}{\#(\{|\pi|\text{ indices}\})} = \lim_{j \to \infty} \frac{N_j |\pi|!}{|\tau[\gamma_j]|^{|\pi|}},
\]
where $N_j$ is the number of $|\pi|$-point subsets of $\tau[\gamma_j]$ that are order-isomorphic to $\pi$. Note that each $|\pi|$-point subset of $\tau[\gamma_j]$ can be divided among block boundaries in the tensor product, which induces some block decomposition of $\pi$.

We can count $N_j$ by taking a sum of the number of $\pi$-point subsets inducing each block decomposition $(\sigma, b) \in B(\pi)$. Note that there are $\binom{|\tau|}{|\sigma|} \cdot t(\sigma, \tau)$ subsets of indices in $\tau$ that are order-isomorphic to $\sigma$. Each of these is an order-respecting way of mapping each index of $\sigma$ to a block of $\tau[\gamma_j]$. Recall though that each block of $\tau[\gamma_j]$ is order-isomorphic to $\gamma_j$, so the number of $|\pi|$-point subsets of $\tau[\gamma_j]$ that are order-isomorphic to $\pi$ is \[
    N_j = \sum_{(\sigma, b) \in B(\pi)}  \pa*{\binom{|\tau|}{|\sigma|}t(\sigma, \tau) \cdot \prod_{i=1}^{|\sigma|} \binom{|\gamma_j|}{|b_i|} t(b_i, \gamma_j) }.
\]
Substituting this in to the previous expression yields \begin{align*}
    \lim_{j \to \infty} t(\pi, \tau[\gamma_j])
    &= \lim_{j \to \infty} \frac{|\pi|!}{|\tau[\gamma_j]|^{|\pi|}} \sum_{(\sigma, b) \in B(\pi)}  \pa*{\binom{|\tau|}{|\sigma|}t(\sigma, \tau) \cdot \prod_{i=1}^{|\sigma|} \binom{|\gamma_j|}{|b_i|} t(b_i, \gamma_j) } \\
    &= \frac{|\pi|!}{|\tau|^{|\pi|}} \sum_{(\sigma, b) \in B(\pi)}  \pa*{\binom{|\tau|}{|\sigma|}t(\sigma, \tau) \cdot \lim_{j \to \infty} \prod_{i=1}^{|\sigma|} \frac{\binom{|\gamma_j|}{|b_i|}}{|\gamma_j|^{|\pi|}} t(b_i, \gamma_j) } \\
    &= \frac{|\pi|!}{|\tau|^{|\pi|}} \sum_{(\sigma, b) \in B(\pi)}  \pa*{\binom{|\tau|}{|\sigma|}t(\sigma, \tau) \cdot \prod_{i=1}^{|\sigma|} \frac{\lim_{j \to \infty} t(b_i, \gamma_j)}{|b_i|!} }.
\end{align*}
This also shows that the limit density converges, as $\lim_{j \to \infty} t(b_i, \gamma_j)$ exists.
\end{proof}

\section{Limit Densities in Random Inflations}

\label{sec:density-inflation}
A special case of the previous theorem is of interest to us, where we take the sequence of tensor products of a fixed permutation $\tau$ with random permutations $\lambda_1, \lambda_2, \ldots$ of increasing lengths. For each $\ell$, $\lambda_\ell$ is sampled from a uniform distribution over all permutations of length $\ell$. We call this sequence of random tensor products the \emph{random inflation} of $\tau$.

\begin{definition}[Random inflation]
The \emph{random inflation} of $\tau$, denoted by $\tau[\bullet]$, is a sequence of random permutations $\tau[\lambda_1], \tau[\lambda_2], \ldots$, where for all $\ell \geq 1$, $\lambda_\ell$ is selected from the uniform distribution over $S_\ell$.\footnote{Unlike the previous definitions of \emph{inflation} and \emph{tensor product}, which produce a single finite permutation, the \emph{random inflation} operation produces a convergent sequence of random permutations.}
\end{definition}

Note that for any pattern $\pi$, by symmetry $\mathbf E[t(\pi, \lambda_j)] = \frac 1{|\pi|!}$ whenever $j \geq |\pi|$. We can actually show something stronger: the sequence $\lambda_1, \lambda_2, \ldots$ is almost surely convergent. This follows from standard concentration inequalities, as shown below.

\begin{lemma}
For any $\pi \in S_k$ and $n \geq k$, let $\lambda_n$ be a random permutation chosen from the uniform distribution on $S_n$. Then, \[
    \Pr\pa*{\left|t(\pi, \lambda_n) - \frac 1{k!}\right| \geq \epsilon} \leq 2e^{-2\epsilon^2 n/k^2}.
\]
\end{lemma}
\begin{proof}
Suppose that we generate $\lambda_n$ by the following random process: choose $n$ independent, identically distributed variables $X_1, \ldots, X_n \sim \operatorname{Uniform}(0, 1)$. With probability $1$, these values are all distinct, so let $\lambda_n$ be the unique permutation that is order-isomorphic to $(X_1, \ldots, X_n)$. We can also define the function $f(X_1, \ldots, X_n) = t(\pi, \lambda_n)$.

Note that for any $i$ and distinct values $x_1, \ldots, x_i, x_i', \ldots x_n$, we have that \[
    |f(x_1, \ldots, x_i, \ldots, x_n) - f(x_1, \ldots, x_i', \ldots, x_n)| \leq \frac{\binom{n-1}{k-1}}{\binom nk} = \frac{k}{n}.
\]
This is because changing $x_i$ can only affect the induced subpermutations of $\lambda_n$ that include the value $\lambda_n(i)$. Then, by McDiarmid's inequality, for any $\epsilon > 0$ \[
    \Pr(|f(X_1, \ldots, X_n) - \mathbf E[f(X_1, \ldots, X_n)]| \geq \epsilon) \leq 2e^{-\frac{2\epsilon^2}{n(k/n)^2}} = 2e^{-2\epsilon^2 n/k^2}.
\]
The result follows from the definition of $f$.
\end{proof}

Since the probability of $t(\pi, \lambda_n)$ deviating from its expected value by more than $\epsilon$ decreases exponentially with $n$, it follows this sequence converges with probability $1$, i.e., \[
    \lim_{j \to \infty} t(\pi, \lambda_j) = \frac 1{|\pi|!}.
\]
We can then apply \cref{thm:tensor-limit-density} to the case of random inflation by plugging in $t(\alpha, \gamma) = \frac{1}{|\alpha|!}$, which yields the following result.

\begin{theorem}
\label{thm:inflation-limit-density}
If $\{ \lambda_j \}$ is a sequence of random permutations with lengths tending to infinity, then the limit density of $\pi$ in the sequence of permutations $\tau[\lambda_1], \tau[\lambda_2], \ldots$ is
\[
    t(\pi, \tau[\bullet]) =
    \frac{|\pi|!}{|\tau|^{|\pi|}} \sum_{(\sigma, b) \in B(\pi)} \left(\binom{|\tau|}{|\sigma|} t(\sigma, \tau) \cdot \prod_{i=1}^{|\sigma|} \frac{1}{|b_i|!^2}\right).
\]
\end{theorem}

\subsection{Example of 3-Point Densities in Random Inflation}

As we are primarily interested in the limit densities of $3$-patterns, for our next example we choose $\pi = 132$. Let $\tau$ be some permutation of length $n$, and let $\pi = 132$. The permutation $\pi$ admits 3 block decompositions, so the limit density of $\pi$ in the random inflation of $\tau$ is \begin{align*}
    t(\pi, \tau[\bullet])
    &= \frac{3!}{n^3} \left(
        \binom{n}{1}t(1, \tau) \cdot \frac{1}{3!^2} +
        \binom{n}{2}t(12, \tau) \cdot \frac{1}{1!^2 2!^2} +
        \binom{n}{3}t(132, \tau) \cdot \frac{1}{1!^2 1!^2 1!^2}
    \right) \\
    &= \frac{3!}{n^3}
    \left(
        \frac{1}{36}\binom{n}{1} +
        \frac{1}{4}\binom{n}{2}t(12, \tau) +
        \binom{n}{3}t(132, \tau)
    \right).
\end{align*}
Notice that there exists additional structure if we focus our attention specifically on $3$-patterns. For the sake of brevity, denote $\frac{3!}{n^3} \binom{n}{3}$ by $a(n)$, $\frac{3!}{n^3}\frac{1}{4}\binom{n}{2}$ by $b(n)$, and $\frac{3!}{n^3}\frac 1{36}\binom n1$ by $c(n)$. Substituting yields \[
    t(132,\tau[\bullet]) = a(n)t(132,\tau) + b(n) t(12,\tau)+ c(n).
\]
We can compute the limit densities of $213$, $312$ and $231$ in the random inflation of $\tau$ in a similar manner, and we see that: \begin{align*}
    t(213,\tau[\bullet]) &= a(n)t(213,\tau) + b(n) t(12,\tau) + c(n), \\
    t(312,\tau[\bullet]) &= a(n)t(312,\tau) + b(n) t(21,\tau) + c(n), \\
    t(231,\tau[\bullet]) &= a(n)t(231,\tau) + b(n) t(21,\tau) + c(n).
\end{align*}
The reason for this similarity is each of the permutations $213$, $312$, and $132$ all have a similar block decomposition structure to $231$. Each has three distinct decompositions: one block of size $3$, two blocks of sizes $1$ and $2$, and three blocks of size $1$.

Similarly, permutations 123 and 321 follow similar formulas, except the coefficient of the second term is doubled since there are two ways to partition each of $123$ and $321$ into two blocks. These limit densities are given by \begin{align*}
    t(123,\tau[\bullet]) &= a(n)t(123,\tau) + 2 b(n) t(12,\tau) + c(n), \\
    t(321,\tau[\bullet]) &= a(n)t(321,\tau) + 2 b(n) t(12,\tau) + c(n).
\end{align*}
To see an example of this symmetry, consider the length-$9$ permutation $\tau = 472951836$. The density of $132$ in $\tau$ is $17/84$, and the density of $12$ in $\tau$ is $1/2$. Therefore, \[
    t(132, \tau[\bullet]) = \frac{29}{162}.
\]
Similarly, we can compute \begin{gather*}
    t(213, \tau[\bullet]) = t(231, \tau[\bullet]) = t(312, \tau[\bullet]) =  \frac{29}{162}, \\
    t(123, \tau[\bullet]) = t(321, \tau[\bullet]) = \frac{23}{162}.
\end{gather*}

\section[k-Inflatability and k-Symmetry]{$k$-Inflatability and $k$-Symmetry}
\label{sec:ksymmetry}

We now turn our attention to inflatable permutation sequences, which are related to quasirandom permutation sequences. First, recall the definitions from \cite{cooper08}.

\begin{definition}[\textit{k}-symmetric]
A convergent permutation sequence $\{\tau_j\}$ is called \textit{$k$-symmetric} if for every $k$-pattern $\pi$, the limit density of $\pi$ in $\{\tau_j\}$ is $1/k!$.
\end{definition}

Note that if a permutation sequence is quasirandom (\cref{def:quasirandom}), then it is $k$-symmetric for all values of $k$. Also, as shown in \cite{kral12}, any $4$-symmetric permutation sequence is quasirandom.

\begin{definition}[\textit{k}-inflatable]
A permutation $\tau$ is called \textit{$k$-inflatable} if the random inflation of $\tau$ is $k$-symmetric.
\end{definition}

As a direct corollary of \cref{thm:inflation-limit-density}, a permutation $\tau$ is $2$-inflatable if and only if $t(12, \tau) = 1/2$, so checking for $2$-inflatability simply amounts to inversion counting. Also, since a permutation sequence is quasirandom if and only if it is $4$-symmetric, there do not exist any nontrivial $4$-inflatable permutations. The only interesting case is $k = 3$.

It was claimed in \cite{cooper08} that the smallest $3$-inflatable permutations are of length $9$, and there are four examples: $472951836$, its inverse, and reflections. However, this permutation is actually the same as $\tau$ from the last section, where we showed by calculation that $472951836$ is \textbf{not} $3$-inflatable. We independently verified our density calculations by conducting experiments where we inflated this permutation with random permutations in $S_{30}$. The code can be found in \cprotect{\href{https://github.com/ekzhang/inflatable/blob/master/random_simulation.py}}{\verb|random_simulation.py|}.

\subsection[Characterization of 3-Inflatable Permutations]{Characterization of $3$-Inflatable Permutations}

Applying \cref{thm:inflation-limit-density} allows us to completely characterize $3$-inflatable permutations.

\begin{theorem}
A permutation $\tau$ of length $n$ is $3$-inflatable if and only if \begin{gather*}
t(12, \tau) = 1/2, \\
t(123, \tau) = t(321, \tau) = (2 n - 7)/(12 (n - 2)), \\
t(132, \tau) = t(213, \tau) = t(231, \tau) = t(312, \tau) = (4 n - 5)/(24 (n - 2)).
\end{gather*}
\end{theorem}

\begin{proof}
First, a 3-inflatable permutation must also be 2-inflatable, so the density of $12$ in $\tau$ must be $1/2$. A specific application of \cref{thm:inflation-limit-density} for three-point densities yields \[
    t(123, \tau[\bullet]) = \frac{9 \cdot 8 \cdot 7}{9^3} t(123, \tau) +  \binom 92 \cdot \frac{3}{9^3} \cdot 2 \cdot t(12, \tau) \cdot \frac 12 + \frac{9}{9^3} \frac{1}{3!}.
\]
The density of inversions should be $\frac 12$, so this is equivalent to \[
    \frac{1}{6} = \frac{n(n-1)(n-2)}{n^3} t(123, \tau) +  \binom n2 \cdot \frac{3}{n^3} \cdot 2 \cdot 
    \frac 12 \cdot \frac 12 + \frac{n}{n^3} \frac{1}{3!}.
\]
Multiplying by $12n^2$ and rearranging yields \[
    2n^2 - 9n + 7 = 12(n-1)(n-2) t(123, \tau).
\]
This means that we must have \[
    t(123, \tau) = \frac{2n^2 - 9n + 7}{12(n-1)(n-2)} = \frac{(2n-7)(n-1)}{12(n-1)(n-2)} = \frac{2n-7}{12(n-2)}.
\]
We can analogously derive the densities of the five other $3$-patterns.
\end{proof}

For some values of $n$, due to divisibility reasons, it may not be possible for any permutation $\tau$ of length $n$ to satisfy the conditions required to be $3$-inflatable. This gives us a way to further narrow down the possible lengths of $3$-inflatable permutations.

\begin{corollary}
If a permutation of length $n$ is $3$-inflatable, then:
\begin{enumerate}
\item $\binom n2$ is even, and
\item $\binom n3$ is divisible by the reduced denominators of $\frac{2n-7}{12(n-2)}$ and $\frac{4 n - 5}{24 (n - 2)}$.
\end{enumerate}
\end{corollary}

This can be described more explicitly by computing the complete set of admissible residue classes modulo $144$, for the lengths of a $3$-inflatable permutation.

\begin{lemma}
For any 3-inflatable permutation $\tau$ with length $n$, \[
    n \equiv 0, 1, 17, 64, 80, 81 \pmod{144}.
\]
\end{lemma}

\begin{proof}
Let $k$ be the number of occurrences of $132$ in the permutation $\tau$. By the above criteria, we have that \[
    \frac{k}{\binom{n}{3}} = \frac{4 n - 5}{24 (n - 2)},
\] \[
    144 k = n(n-1)(4 n - 5),
\] \[
    144 \mid n(n-1)(4n-5).
\]
Also, \[
    \frac{k}{\binom{n}{3}} = \frac{2 n - 7}{12 (n - 2)},
\] \[
    72k = n(n - 1)(2 n - 7),
\] \[
    72 \mid n(n - 1)(2 n - 7).
\]
Together, these two divisibility requirements imply the result.
\end{proof}

Therefore, the smallest possible nontrivial length for a $3$-inflatable permutation is $17$. In this case, $\binom{17}{2} = 136$ and $\binom{17}{3} = 680$, while $(2n-7)/(12(n-2)) = 3/20$ and $(4 n - 5)/(24 (n - 2)) = 7/40$. This means that in a $3$-inflatable permutation of length $17$, the number of occurrences of the patterns $123$ and $321$ must both be $102$, while the four other $3$-patterns should each occur $119$ times.

\begin{figure}
    \centering
    \begin{asy}
    unitsize(5cm);
    drawperm("E534BGA9HC2D1687F");
    label("E534BGA9HC2D1687F", (0.5, -0.03), S);
    drawperm("G54ABC319HF678ED2", 1.2);
    label("G54ABC319HF678ED2", (0.5 + 1.2, -0.03), S);
    \end{asy}
    \caption{Plots of selected minimal-length 3-inflatable permutations.}
    \label{fig:plot-length17}
\end{figure}

Examples of these minimum-length 3-inflatable permutations were found through an optimized computer search. For example, the permutations \[\text{E534BGA9HC2D1687F and  G54ABC319HF678ED2}\] are both 3-inflatable, where capital letters denote indices greater than nine. Visual plots of these permutations are shown in \cref{fig:plot-length17}.

\subsection{Centrally Symmetric Permutations}

Now we want to define a property of certain permutations that is useful for computer search, as it increases the likelihood of being symmetric. Let us define the \emph{rotation} of a permutation $\pi$, denoted $R(\pi)$, to be a permutation of the same length such that $R(\pi)(i) = n + 1 - \pi(n + 1 - i)$ for all $i$. This can be visualized by drawing the permutation on a square grid, then rotating it by 180 degrees about its center.

We call a permutation $\pi$ of length $n$ \emph{centrally symmetric} if it is equal to its rotation $R(\pi)$. In other words, $\pi(i) + \pi(j) = n+1$ whenever $i+j = n + 1$. Examples of centrally symmetric permuations are shown in \cref{fig:centrally-symmetric}, as well as the right plot in \cref{fig:plot-length17}. The importance of centrally symmetric permutations is explained by the following lemma.

\begin{figure}
    \centering
    \begin{asy}
    unitsize(3.5cm);
    drawperm("123456");
    label("123456", (0.5, -0.03), S);
    drawperm("924753681", 1.2);
    label("924753681", (0.5 + 1.2, -0.03), S);
    drawperm("482156A937", 2.4);
    label("482156A937", (0.5 + 2.4, -0.03), S);
    \end{asy}
    \caption{Examples of small centrally symmetric permutations.}
    \label{fig:centrally-symmetric}
\end{figure}

\begin{lemma}
For permutations $\pi$ and $\gamma$, \[
    t(\pi,\gamma) = t(R(\pi),R(\gamma)).
\]
\end{lemma}

\begin{proof}
Let $k = |\pi|$, and let $T_1$ be the set of $k$-point subsets of $\gamma$ that are order-isomorphic to $\pi$, while $T_2$ is the set of $k$-point subsets of $R(\gamma)$ that are order-isomorphic to $R(\pi)$. For any set of indices $\{ i_1, i_2, \dots, i_k \}$ in $\gamma$, there is a corresponding set of indices, $\{ n + 1 - i_1, n + 1 - i_2, \dots, n + 1 - i_k \}$ in $R(\gamma)$. If the former indices induce some permutation $\pi$ in $\gamma$, then the latter set induces $R(\pi)$ in $R(\gamma)$ from the definition of rotation. This gives a bijection between $T_1$ and $T_2$, so the result follows.
\end{proof}

\begin{corollary}
If $\pi$ is centrally symmetric, then
\[t(\pi,\gamma) = t(\pi,R(\gamma)).\]
\end{corollary}

That means centrally symmetric permutations automatically give us some equalities between densities. For example, if $\pi$ is centrally symmetric, then \[
    t(132,\pi) = t(231,\pi) \quad \text{and} \quad t(312,\pi) = t(213,\pi),
\]
which explains why it was easier to find examples of $3$-inflatable permutations when we restricted our search to random centrally symmetric permutations.

Although it was intractable to check all $17!$ permutations for the $3$-inflatability criterion, we could reduce the search space by checking only those permutations $\pi$ that were centrally symmetric. Empirically, we found that adding this constraint greatly increased the proportion of hits. Of the $8 \cdot 2^8$ centrally symmetric permutations of length $17$, a computer search found that exactly $750$ of them are $3$-inflatable. A list of these $750$ permutations is given in the file \cprotect{\href{https://github.com/ekzhang/inflatable/blob/master/search/inflatable_list17.txt}}{\verb|inflatable_list17.txt|}.

Furthermore, additional randomized computer searches found many more examples of centrally symmetric, $3$-inflatable permutations of lengths $64$, $80$, $81$, $144$, and $145$. These lists are non-exhaustive, and the computer search can easily be extended to larger lengths, but we stopped at this point due to complexity. For example, from the results in the \cprotect{\href{https://github.com/ekzhang/inflatable/tree/master/search}}{\verb|search/|} directory, we estimate there to be on the order of $10^{37}$ centrally symmetric $3$-inflatable permutations of length $64$, and $10^{115}$ centrally symmetric $3$-inflatable permutations of length $145$.

Regardless, this covers all admissible residue classes of lengths modulo $144$, and there seems to be a large number of such permutations of each length. We would like to suggest the following open conjecture for further study.

\begin{conjecture}
For any positive integer $x$ belonging to an admissible residue class of lengths modulo $144$, there exists a $3$-inflatable permutation of length $x$.
\end{conjecture}

\section{Structure of 3-Inflatable Permutations}
\label{sec:structure3}

To explain some of the empirical results of our computer search, we found an interesting structural property of $3$-inflatable permutations.

\begin{theorem}
If $\tau_1$ and $\tau_2$ are both $k$-inflatable, then $\tau_1[\tau_2]$ is also $k$-inflatable.
\end{theorem}
\begin{proof}
Note that tensor multiplication is associative, so for any permutation $\gamma$, we must have $(\tau_1[\tau_2])[\gamma] = \tau_1[\tau_2[\gamma]]$. Also, for any length-$k$ permutation $\pi$, \[
    \lim_{|\gamma| \to \infty} t(\pi, \tau_2[\gamma]) = \frac{1}{k!}.
\]
Since our current discussion only considers densities of pattern permutations with length at most $k$, we can substitute $\tau_2[\gamma]$ for a random permutation $\gamma_2$ of the same length without changing the expected value of the expression. Thus, \[
    \lim_{|\gamma| \to \infty} t(\pi, \tau_1[\tau_2[\gamma]]) = \lim_{|\gamma_2| \to \infty} t(\pi, \tau_1[\gamma_2]) = \frac{1}{k!},
\]
so $\tau_1[\tau_2]$ is $k$-inflatable.
\end{proof}

This result can be used to generate an infinite number of $3$-inflatable permutations with increasing lengths. For example, we can use any $3$-inflatable permutations of size $17$ to construct $3$-inflatable permutations of size $17^n$ for any $n$.

It is also natural to check if the set of admissible lengths of $3$-inflatable permutations (modulo $144$) is closed under multiplication. This is the case, as shown in \cref{fig:mult_table}. We are currently unaware of other operations that preserve $k$-inflatability.

\begin{figure}
    \centering

    \begin{tabular}{|c|cccccc|}
    \hline
     $\times$  & 0 &  1 & 17 & 64 & 80 & 81 \\ \hline
     0  & 0 &  0 &  0 &  0 &  0 &  0 \\
     1  & 0 &  1 & 17 & 64 & 80 & 81 \\
     17 & 0 & 17 &  1 & 80 & 64 & 81 \\
     64 & 0 & 64 & 80 & 64 & 80 &  0 \\
     80 & 0 & 80 & 64 & 80 & 64 &  0 \\
     81 & 0 & 81 & 81 &  0 &  0 & 81 \\
    \hline
    \end{tabular}

    \caption{Multiplication table of lengths of $3$-inflatable permutations, modulo $144$.}
    \label{fig:mult_table}
\end{figure}

\section{Acknowledgements}

We would like to thank Yufei Zhao, for suggesting the original topic of this project; the anonymous reviewer, for their thorough review and suggestions that we incorporated into this paper; and
the PRIMES program, for giving us the opportunity to do this research.

\bibliographystyle{alpha}
\bibliography{inflatable}

\end{document}